\documentclass{amsart}
\usepackage{graphicx, color}
\usepackage{amscd}
\usepackage{amsmath}
\usepackage{amsfonts}
\usepackage{amssymb}
\usepackage{mathrsfs}

\usepackage{srcltx}
\usepackage{enumitem}
\usepackage{amssymb, color}
\usepackage{xcolor}
\usepackage{hyperref}
\usepackage{amsmath}
\usepackage{pifont}
\usepackage[utf8]{inputenc}
\newtheorem{theorem}{Theorem}[section]
\newtheorem{lemma}[theorem]{Lemma}
\newtheorem{prop}[theorem]{Proposition}
\newtheorem{rem}[theorem]{Remark}

\newtheorem{definition}{Definition} [section]

\allowdisplaybreaks
\numberwithin{equation}{section}

\numberwithin{equation}{section} 

\renewcommand{\leq}{\leqslant}

\renewcommand{\geq}{\geqslant}
\newcommand{\R}{\mathbb R}

\def \O{\Omega}  
\begin{document} 
\title[Systems arising in electromagnetism]{Existence of solutions for systems \\ arising in electromagnetism} 

\author[M.K. Hamdani]{M.K. Hamdani}
\address[M.K. Hamdani] {Military School of Aeronautical Specialities, Sfax \&
 Mathematics Department, Faculty of Science, University of Sfax, Sfax,  Tunisia.} \email{hamdanikarim42@gmail.com} 

\author[D.D. Repov\v{s}]{ D.D. Repov\v{s}}
\address[D.D. Repov\v{s}] {Faculty of Education and Faculty of Mathematics and Physics, University of Ljubljana \& Institute of Mathematics, Physics and Mechanics,  Ljubljana, Slovenia}
\email{dusan.repovs@guest.arnes.si}
 
\begin{abstract}
In this paper, we study the following $p(x)$-curl systems:
\begin{eqnarray*}
\begin{cases}
\nabla\times(|\nabla\times \mathbf{u}|^{p(x)-2}\nabla\times \mathbf{u})+a(x)|\mathbf{u}|^{p(x)-2}\mathbf{u}=\lambda f(x,\mathbf{u})+\mu g(x,\mathbf{u}),\quad\nabla\cdot \mathbf{u}=0,\; \mbox{ in } \Omega, \\
|\nabla\times \mathbf{u}|^{p(x)-2}\nabla\times \mathbf{u}\times \mathbf{n}=0,\quad \mathbf{u}\cdot \mathbf{n}=0, \mbox{ on } \partial\Omega,
\end{cases}
\end{eqnarray*}
where 
$\Omega\subset\mathbb{R}^3$ 
is a bounded simply connected domain with a
 $C^{1,1}$-boundary, denoted by
  $\partial\Omega$, 
  $p:\overline{\Omega}\to (1, +\infty)$
   is a continuous function,
$a \in L^\infty(\Omega)$,
 $f,g : \Omega \times \mathbb{R}^3\to \mathbb{R}^3$ 
 are Carath\'{e}odory functions, and 
 $\lambda,\mu$ are two parameters. 
 Using variational arguments based on Fountain theorem and Dual Fountain theorem, we establish some existence and non-existence results for solutions of this problem. Our main results generalize the results of Xiang et al. (J. Math. Anal. Appl., 2017), Bahrouni and Repov\v{s}  (Complex Var. Elliptic Equ., 2018), and  Ge and Lu (Mediterr. J. Math., 2019).
\end{abstract}


\maketitle


\section{Introduction}\label{sect1}
Since the variable exponent spaces have been thoroughly
studied by  Kov\'a\u{c}ik-R\'akosn\'ik \cite{KR},
they have been used in  previous decades to model various phenomena.
In the studies of a class of non-standard variational problems and PDEs,
variable exponent spaces play an important role, e.g. in electrorheological fluids \cite{RR1996,RR2001,R}, thermorheological fluids \cite{AR,AS}, and image processing \cite{AMS,CLR,LLP}. For nonlinear problems with variable growth,  there has been a great deal of interest in studying the existence, multiplicity, uniqueness and regularity of solutions - for the main results (as well as definitions of some these properties) see
\cite{ACX,AH2,AMC,BO,Chung1,Chung2,NTC1,JF,FZ,FanZ,GL,Ham,KR,MR,MR2007,Dusan,RD} and the references therein.

The $p(x)$-curl operator defined by $\nabla\times(|\nabla\times \mathbf{u}|^{p(x)-2}\nabla\times \mathbf{u})$ is a generalization of the $p$-curl operator
  in which the constant exponent $p$ has been
replaced by a variable exponent $p(x)$. The $p(x)$-curl systems possess more complicated structure than the $p$-curl operators, due to the fact that they are not homogeneous.  Therefore the study of various mathematical problems with variable exponent is very interesting and raises many difficult mathematical problems.

Moreover, the study of nonlinear elliptic equations involving quasilinear homogeneous type operators like the $p$-Laplace or $p$-curl operators is based on the theory of standard Sobolev spaces $W^{1,p}(\O)$ in order to find weak solutions - see \cite {AMS0,H2,Harrabi,HM}. These spaces consist of functions that have weak derivatives and satisfy certain integrability conditions. In the case of nonhomogeneous $p(x)$-Laplace operators, the natural setting for this approach is to use of the variable exponent Sobolev spaces. The basic idea is to replace the Lebesgue spaces $L^p(\O)$ by more general spaces $L^{p(\cdot)}(\O)$, called the variable exponent Lebesgue spaces.  However, in literature the only results involving the $p(x)$-curl systems by variational methods can be found in \cite{AMS,BD,BJF,XWZ}.

Throughout the paper, vector functions and spaces of vector functions will be in boldface. We shall use $\partial_x$ to denote the partial derivative of a function with respect to the variable $x$.

To introduce our problem precisely, we first give some notations. Let $\mathbf{u} = (u_1, u_2, u_3)$ be a vector
function on $\O$. The divergence of $\mathbf{u}$ is denoted by
$$\nabla\cdot \mathbf{u}=\partial_{x_1}u_1+\partial_{x_2}u_2+\partial_{x_3}u_3$$
and the curl of $\mathbf{u}$, written
$curl\;\mathbf{u}$ or $\nabla\times \mathbf{u}$, is defined to be the vector field
$$ \nabla\times \mathbf{u}=\left\langle \partial_{x_2}u_3- \partial_{x_3}u_2,\partial_{x_3}u_1-\partial_{x_1}u_3, \partial_{x_1}u_2 - \partial_{x_2}u_1\right\rangle.$$
Throughout this paper, unless otherwise stated, we shall always assume that exponent $p(x)$ is continuous
on $\overline{\O}$ with
$$1<p^-=\min_{x\in\overline{\O}}p(x)\leq p^+=\max_{x\in\overline{\O}}p(x)<3,$$
and satisfies the logarithmic continuity, i.e. that there exists a function $\omega : \R_0^+\to\R_0^+$ such that
for all $ \  x,y\in \overline{\O},$
\begin{eqnarray}\label{cond12}
 x-y|<1,\;|p(x)-p(y)|\leq \omega(|x-y|), \mbox{ and } \lim_{\tau\to 0^+}\omega(\tau)\log\frac{1}{\tau}=C<\infty.
\end{eqnarray}
In 2016,  Antontsev-Miranda-Santos \cite{Antontsev} studied the qualitative properties of solutions for the following $p(x,t)$-curl systems:
\begin{eqnarray} \label{Problem2}
\begin{cases}
\partial_t \mathbf{u}+\nabla\times(|\nabla\times \mathbf{u}|^{p(x,t)-2}\nabla\times \mathbf{u})=f(\mathbf{u}),\quad\nabla\cdot \mathbf{u}=0,\mbox{ in } \Omega\times (0,T), \\
|\nabla\times \mathbf{u}|^{p(x,t)-2}\nabla\times \mathbf{u}\times \mathbf{n}=0,\quad \mathbf{u}\cdot \mathbf{n}=0, \mbox{ on } \partial\Omega\times (0,T),\\
\mathbf{u}(x,0)=\mathbf{u}_0(x),\mbox{ in }\Omega,\\
\end{cases}
\end{eqnarray}
where $\nabla\times(|\nabla\times \mathbf{u}|^{p(x,t)-2}\nabla\times \mathbf{u})$  is the $p(x,t)$-curl operator $$f(\mathbf{u})=\lambda \mathbf{u} (\int_\O|\mathbf{u}|^2dx)^\frac{\rho-2}{2} \ \  \mbox{where} \ \
\lambda\in \{-1, 0, 1\} \ \ \mbox{and} \ \   \rho  \mbox{ is a positive constant}.$$ They introduced a suitable functional framework and a
convenient basis  in order to apply  Galerkin's method and they studied the blow-up and finite time extinction properties of solutions, depending on the values of $\lambda$ and $\rho$. In the same year, Xiang-Wang-Zhang \cite{XWZ} used for the first time, the variational methods for equations involving $p(x)$-curl operator of the following type:
\begin{eqnarray} \label{Xiang}
\begin{cases}
\nabla\times(|\nabla\times \mathbf{u}|^{p(x)-2}\nabla\times \mathbf{u})+a(x)|\mathbf{u}|^{p(x)-2}\mathbf{u}=f(x,\mathbf{u}),\quad \nabla\cdot \mathbf{u}=0, \mbox{ in } \Omega, \\
|\nabla\times \mathbf{u}|^{p(x)-2}\nabla\times \mathbf{u}\times \mathbf{n}=0, \quad\mathbf{u}\cdot \mathbf{n}=0, \mbox{ on } \partial\Omega.
\end{cases}
\end{eqnarray}
They studied  the existence and multiplicity of solutions for system \eqref{Xiang} with the following assumptions on $a(x)$ and $f(x,\mathbf{u})$:
\begin{enumerate}
  \item[$\mathcal{(A)}:$] $a(x)\in L^\infty(\O)$ and
   there exist $a_0,a_1>0$ such that $a_0\leq a(x)\leq a_1$  for all $x\in \O$.
\item[$(H_1):$] There exists $F : \O\times \R^3 \to \R$ which is differentiable with respect to $\mathbf{u}\in \R^3$ and such that $$f(x,\mathbf{u})= \partial_\mathbf{u} F(x,\mathbf{u}): \O\times \R^3 \to \R^3$$ is a Carath\'{e}odory function.
\item[$(H_2):$] There exist $C> 0$,
 $q\in C(\overline{\Omega})$,
  and $$1 <  q(x)< p^*(x)=\frac{3p(x)}{3-p(x)} \ \  \mbox{in} \ \ \overline{\Omega}$$ such that
$$ |f(x,\mathbf{u})|\leq C(1+|\mathbf{u}|^{q(x)-1}),~~ \mbox{ for all } (x,\mathbf{u}) \in {\Omega}\times {\R^3}.$$
\item[$(H_3):$] There exists a constant $\mu > p^+$ such that $$0 < \mu F(x,\mathbf{u}) \leq f(x,\mathbf{u})\cdot\mathbf{u}~~ \mbox{ for all } x \in \Omega \mbox{ and } \mathbf{u}\in \R^3\setminus\{0\}.$$
\item[$(H_4):$] $\limsup_{\mathbf{u}\rightarrow 0}\frac{|f(x,\mathbf{u})|}{|\mathbf{u}|^{p(x)-1}}=0 \mbox{ uniformly in } x \in \O.$
\item[$(H_5):$] $\inf_{x\in\O,\mathbf{u}\in\R^3, |\mathbf{u}|=1}F(x,\mathbf{u})>0.$
\item [$(H_6):$] $F(x,-\mathbf{u})=F(x,\mathbf{u})$ \;\;for all $(x,\mathbf{u})\in\O\times\R^3.$
\end{enumerate}
The proofs in \cite{XWZ} are based on Mountain Pass theorem and Symmetric Mountain Pass theorem. Under the conditions $\mathcal{(A)}$ and $(H_1)-(H_6)$, the following was proved in \cite{XWZ}.
\begin{theorem}(see \cite[Theorems $1.1$ and $1.2$]{XWZ}). Suppose that  $$p(x) <
q(x) <\frac{3p(x)}{3-p(x)} \ \ \mbox{ for all} \ \  x\in \overline{\O}.$$ Then the following holds:
\begin{enumerate}
  \item If $a(x)$ satisfies $\mathcal{(A)}$ and $f(x,\mathbf{u})$ satisfies $(H_1)-(H_5)$,
   then system \eqref{Xiang} has one nontrivial mountain pass solution.
\item If $a(x)$ satisfies $\mathcal{(A)}$ and $f(x,\mathbf{u})$ satisfies $(H_1)-(H_4)$ and $(H_6)$, then system \eqref{Xiang} has infinitely many nontrivial
mountain pass
solutions.
\end{enumerate}
\end{theorem}
In 2019, with the same  method as in
\cite{XWZ}, Ge-Lu \cite{BJF} gave some weaker conditions than in \cite{XWZ} and they proved the existence and the multiplicity of solutions for \eqref{Xiang}. In $2017$, Bahrouni-Repov\v{s} \cite{BD} studied the following $p(x)$-curl system
\begin{eqnarray} \label{Dusan}
\begin{cases}
\nabla\times(|\nabla\times \mathbf{u}|^{p(x)-2}\nabla\times \mathbf{u})=\lambda f(x,\mathbf{u})-\mu g(x,\mathbf{u}),\quad \nabla\cdot \mathbf{u}=0, \;\mbox{ in } \Omega, \\
|\nabla\times \mathbf{u}|^{p(x)-2}\nabla\times \mathbf{u}\times \mathbf{n}=0,\quad \mathbf{u}\cdot \mathbf{n}=0, \mbox{ on } \partial\Omega.
\end{cases}
\end{eqnarray}
Clearly, this problem is a special case of our main system
when $a\equiv 0$. Bahrouni-Repov\v{s} studied the existence of solutions
for system \eqref{Dusan} when $f$ satisfies $(H_1)$, plus the following conditions:
\begin{enumerate}
\item[$(F_2):$] There exist $\alpha,\beta> 0$ and $q\in C(\overline{\Omega})$ such that $$p^+<q(x)< p^*(x)=\frac{3p(x)}{3-p(x)} \ \ \mbox{ in} \ \  \overline{\Omega}$$ and
$$ |F(x,\mathbf{u})|\geq \alpha |\mathbf{u}|^{q(x)}\mbox{ and } |f(x,\mathbf{u})|\leq\beta(1+|\mathbf{u}|^{q(x)-1}),\;\mbox{ for all } (x,\mathbf{u}) \in \overline{\Omega}\times {\R^3}.$$
\end{enumerate}
whereas for $g$ they made the following assumptions:\begin{enumerate}
\item[$(G_1):$] There exist a nonnegative function $g\in L^\infty(\O)$ and $r\in C(\overline{\Omega})$ such that
    $$p^+<r^-\leq r(x)<q^-\;\mbox{ and }\; G(x,\mathbf{u})=g(x)|\mathbf{u}|^{r(x)}, \;\ \ \mbox{for all} \ \  (x,\mathbf{u}) \in \overline{\Omega}\times {\R^3}.$$
\item[$(G_2):$] $G : \O\times \R^3 \to \R$ is differentiable with respect to $\mathbf{u}\in \R^3$ and $g = \partial_\mathbf{u} G(x,\mathbf{u}):\O\times \R^3 \to \R^3$ is a Carath\'{e}odory function.
\item[$(G_3):$] There exist $\gamma,\mu>0$, $L>1$ and $k,r\in C(\overline{\Omega})$ such that $1 < k < p^-$ and
$1 < r(x) < p^*(x)$,
$$|g(x,\mathbf{u})|\leq \mu(1+|\mathbf{u}|^{r(x)-1}),\;\;\mbox{ for all } (x,\mathbf{u}) \in \overline{\Omega}\times {\R^3},$$
$$\limsup_{\mathbf{u}\to0}\frac{G(x,\mathbf{u})}{|\mathbf{u}|^{p^+}}=0 \mbox{ uniformly in } x \in \O,$$
and
$$\sup_{\mathbf{u}\in \mathbf{W}^{p(x)}(\O)}\int_\O G(x,\mathbf{u})dx>0,\;|G(x,\mathbf{u})|\leq \gamma|\mathbf{u}|^{k(x)},\mbox{ for all } x\in \O \mbox{ and } |\mathbf{u}|>L. $$
\end{enumerate}
Under the conditions $(H_1),(F_2), (G_1)-(G_3)$, the following was proved in \cite{BD}.
\begin{theorem}
(see \cite[Theorems $1.1$ and $1.2$]{BD})
\begin{enumerate}
  \item Assume that hypotheses $(H_1),(F_2), (G_1)-(G_2)$ hold. Then there exist $\lambda_1, \mu_1 > 0$ such that, if $0 < \lambda < \lambda_1$ and $\mu > \mu_1$, then system \eqref{Dusan} does
not have any nontrivial weak solutions.
  \item Assume that hypotheses $(H_1),(F_2), (G_1)-(G_2)$ hold. Then for each $\mu > 0$, there exists $\lambda_\mu > 0$ such that if $\lambda>\lambda_\mu$, then system \eqref{Dusan} has at
least one nontrivial weak solution.
\item Assume that hypotheses $(H_1),(F_2), (G_2)-(G_3)$ hold. Then there exist $\lambda_2, \lambda_3, r > 0$
such that, if $\lambda\in[\lambda_2,\lambda_3]$, then there exists $\mu_2 > 0$ with the following property: for each $\mu\in[0,\mu_2]$,
system \eqref{Dusan} has at least three solutions whose norms are less than $r$.
\end{enumerate}
\end{theorem}

Motivated by these results, we study in this paper the existence of solutions for the following $p(x)$-curl systems by means of Fountain theorem and Dual
Fountain theorem.
{\small
\begin{eqnarray} \label{main}
\begin{cases}
\nabla\times(|\nabla\times \mathbf{u}|^{p(x)-2}\nabla\times \mathbf{u})+a(x)|\mathbf{u}|^{p(x)-2}\mathbf{u}=\lambda f(x,\mathbf{u})+\mu g(x,\mathbf{u}),\quad\nabla\cdot \mathbf{u}=0,\; \mbox{ in } \Omega, \\
|\nabla\times \mathbf{u}|^{p(x)-2}\nabla\times \mathbf{u}\times \mathbf{n}=0,\quad \mathbf{u}\cdot \mathbf{n}=0, \mbox{ on } \partial\Omega,
\end{cases}
\end{eqnarray}
}
where $\O\subset\R^3$ is a bounded simply connected domain with a $C^{1,1}$
boundary denoted by $\partial\O$, $\lambda$ and $\mu$ are parameters, $p:\overline{\O}\to (1, +\infty)$ is a continuous function,
$a \in L^\infty(\O)$, and $f,g : \O \times \R^3\to \R^3$ are Carath\'{e}odory functions.
\begin{rem}
To the best of our knowledge, there are no results concerning curl systems with variable exponent based on
Fountain theorem and Dual Fountain theorem. In this context, the results of our paper can be seen as a generalization of the results above, to the $p(x)$-curl systems arising in electromagnetism.
\end{rem}
We shall impose the following condition on $a(x)$:
\begin{enumerate}[label=$\mathcal{(A)}\sb{\arabic*}$ :, ref=$\mathcal{(A)}\sb{\arabic*}$]
  \item\label{hypot:a} $a(x)\in L^\infty(\O)$  such that $\inf_{x\in \O}a(x)=a^->0$.
  \end{enumerate}
We shall also assume that $f(x,\mathbf{u})$ and $g(x,\mathbf{u})$  satisfy the following global conditions:
\begin{enumerate}[label=($f\sb{\arabic*}$) :, ref=$f\sb{\arabic*}$]
\item\label{hypot:F1} Condition $(H_1)$ stated above.
\item\label{hypot:F2} There exist $c_1> 0$ and $q(x)\in C(\overline{\Omega})$ such that $$1 < p^+<q^-\leq q(x)< p^*(x)=\frac{3p(x)}{3-p(x)} \ \ \mbox{in} \ \  \overline{\Omega}$$ and
$$ |f(x,\mathbf{u})|\leq c_1|\mathbf{u}|^{q(x)-1},~~ \mbox{ for all } (x,\mathbf{u}) \in \overline{\Omega}\times {\R^3}.$$
\item\label{hypot:F3} There are constants $l >0$ and $\theta > p^+$ such that $$0 < \theta F(x,\mathbf{u}) \leq f(x,\mathbf{u})\cdot\mathbf{u}~~ \mbox{ for all } |\mathbf{u}| \geq l  \mbox{ and } x \in \Omega.$$
\item\label{hypot:F5} Condition $(H_6)$ stated above.
\item\label{hypot:F4} $\limsup_{\mathbf{u}\rightarrow \infty}\frac{F(x,\mathbf{u})}{|\mathbf{u}|^{p^-}}=0 \mbox{ uniformly in } x \in \O.$
\item\label{hypot:F6} $f(x,\mathbf{u})\cdot\mathbf{u} > 0$ for all $(x,\mathbf{u})\in \Omega \times \R^3.$
\end{enumerate}
\begin{enumerate}[label=($g\sb{\arabic*}$) :, ref=$g\sb{\arabic*}$]
\item\label{hypot:G1} Condition $(G_2)$ stated above.
\item\label{hypot:G2} There exist $c_2 > 0$ and $\gamma(x) \in C(\overline{\Omega})$ such that
$$1 <\gamma(x)<\gamma^+<p^-< p^*(x)=\frac{3p(x)}{3-p(x)} \ \  \mbox{in} \ \  \overline{\Omega}$$ and
   $$ |g(x,\mathbf{u})|\leq c_2|\mathbf{u}|^{\gamma(x)-1}~~ \mbox{ for all } (x,\mathbf{u}) \in \overline{\Omega}\times {\R^3}.$$
\item\label{hypot:G4} $G(x,-\mathbf{u}) = G(x,\mathbf{u})$\;\; for all $(x, \mathbf{u}) \in \O \times \R^3$.
\item\label{hypot:G3} $\liminf_{\mathbf{u}\rightarrow 0}\frac{G(x,\mathbf{u})}{|\mathbf{u}|^\alpha}\geq0 \mbox{ uniformly in } x \in \O$ with $0<\alpha<p^-.$
\item\label{hypot:G5} $g(x,\mathbf{u})\cdot\mathbf{u} > 0$ for all $(x,\mathbf{u})\in \Omega \times \R^3.$
\end{enumerate}
The variational structure of this problem leads us to introduce the following space
$$\mathbf{W}^{p(x)}(\O)=\{\mathbf{v}\in \mathbf{L}^{p(x)}(\O):\nabla\times \mathbf{v}\in \mathbf{L}^{p(x)}(\O),\nabla\cdot \mathbf{v}=0,\mathbf{v}\cdot \mathbf{n}\mid_{\partial \O}=0\}, $$
see Section \ref{sect2.1} for more details. Let us proceed with setting system \eqref{main} in the variational structure.
A function $\mathbf{u}\in \mathbf{W}^{p(x)}(\O)$ is said to be a weak solution of \eqref{main} if
$$\int_\O|\nabla\times \mathbf{u}|^{p(x)-2}\nabla\times \mathbf{u}\cdot\nabla\times \mathbf{v}dx +\int_\O a(x)|\mathbf{u}|^{p(x)-2}\mathbf{u}\cdot \mathbf{v}dx =\lambda\int_\O f(x,\mathbf{u})\cdot \mathbf{v}dx+\mu\int_\O g(x,\mathbf{u})\cdot \mathbf{v}dx,$$\mbox{for all} $\mathbf{v}\in \mathbf{W}^{p(x)}(\O).$
The Euler-Lagrange functional associated to system \eqref{main} is defined by $$I_{\lambda,\mu}=\Phi-\lambda J-\mu\Psi,\;\;\lambda,\mu\in\R,$$ where
$$\Phi(\mathbf{u})=\int_\O\frac{1}{p(x)}\left( |\nabla\times \mathbf{u}(x)|^{p(x)}+a(x)|\mathbf{u}(x)|^{p(x)}\right)dx,$$ $$J(\mathbf{u})=\int_\O F(x,\mathbf{u})dx,\quad \Psi(\mathbf{u})=\int_\O G(x,\mathbf{u})dx,$$
\\and $$F(x,\mathbf{u})=\int_{0}^{\mathbf{u}}f(x,\mathbf{s})ds,\quad G(x,\mathbf{u})=\int_{0}^{\mathbf{u}}g(x,\mathbf{s})ds.$$

Now, we can state our main results as follows.
\begin{theorem}\label{th1}
Assume that $a(x)$ satisfies \ref{hypot:a}, that conditions
\eqref{hypot:F1}-\eqref{hypot:F5},  \eqref{hypot:G1}-\eqref{hypot:G4} hold, and
that
 $p^+<q^-\leq q(x)\leq p_2^*(x),~~\gamma^+<p^-$.
Then system \eqref{main} has a sequence of weak solutions $(\pm \mathbf{u}_k)$ in $\mathbf{W}^{p(x)}(\O)$ for every $\lambda >0$, ${\mu>0}$, such that  $I_{\lambda,\mu}(\pm \mathbf{u}_k ) \rightarrow +\infty, ~~ as~~ k \rightarrow +\infty$.
  \end{theorem}

\begin{theorem}\label{th2}
Assume that $a(x)$ satisfies \ref{hypot:a}, that conditions
 \eqref{hypot:F1}- \eqref{hypot:F2}, \eqref{hypot:F5}, \eqref{hypot:G1}-\eqref{hypot:G3} hold,
  and that
   $p^+<q^-\leq q(x)\leq p_2^*(x),~~\gamma^+<p^-$.
Then system \eqref{main} has a sequence of weak solutions $(\pm u_k )$ in $\mathbf{W}^{p(x)}(\O)$ for every $\lambda > 0$, $\mu > 0$, such that
  $I_{ \lambda,\mu} (\pm \mathbf{u} _k ) < 0~~ and~~ I_{ \lambda,\mu }(\pm \mathbf{u}_ k ) \rightarrow 0,~~ as~~ k \rightarrow +\infty$.
  \end{theorem}

\begin{theorem}\label{th3}
Assume that $a(x)$ satisfies \ref{hypot:a}, that conditions
 \eqref{hypot:F1}-\eqref{hypot:F2},\eqref{hypot:F4}, \eqref{hypot:G1}-\eqref{hypot:G2}, \eqref{hypot:G3} hold,
  and that
   $p^+<q^-\leq q(x)\leq p_2^*(x),~~\gamma^+<p^-$.
Then system \eqref{main} has at least one nontrivial weak solution in $\mathbf{W}^{p(x)}(\O)$ for every $\lambda < 0,~\mu > 0$.
\end{theorem}
\begin{theorem}\label{th4}
Assume that $a(x)$ satisfies \ref{hypot:a}. If \eqref{hypot:F1}-\eqref{hypot:F2}, \eqref{hypot:F6}, \eqref{hypot:G1}-\eqref{hypot:G2},\eqref{hypot:G5} hold and $p^+<q^-\leq q(x)\leq p_2^*(x),~~\gamma^+<p^-$.
Then system \eqref{main} has no nontrivial weak solution in $\mathbf{W}^{p(x)}(\O)$ for every  $\lambda < 0,~ \mu < 0$.
\end{theorem}

We conclude with an outline of the structure of the paper.
 In Section \ref{sect2}, we introduce some preliminary results
 and  in Section \ref{sect3}, we give the proofs of the main results.
\section{Preliminaries}\label{sect2}
In this section we shall give some preliminary results which will be used in the sequel.
\subsection{\bf Variable exponent Lebesgue and Sobolev spaces}\label{sect2.1}
To study our problems, we shall need to introduce certain function spaces. Denote
$$
C_+(\Omega)=\big\{p\in C(\overline{\Omega}):\min_{x\in\Omega}p(x)>1 \big\}.
$$

 \begin{definition} \label{def2.1} \rm
 The variable exponent Lebesgue space $L^{p(x)}(\Omega)$ is defined by
 $$
L^{p(x)}(\Omega)=\big\{\mathbf{u}; \;\mathbf{u} \mbox{ is a measurable real-valued function such that } \int_\Omega | \mathbf{u} |^{p(x)}dx <+\infty \big\},
$$
and is endowed with the so-called Luxemburg norm
 $$
|\mathbf{u} |_{p(x)}=\inf \big\{ \lambda>0: \int_\Omega \left|\frac{\mathbf{u}}
 {\lambda}\right|^{p(x)}dx \leq 1 \big\}.
$$
 \end{definition}
If $p(x)=p\equiv$
constant for every $x\in \Omega$, then the
$L^{p(x)}(\Omega)$ space is reduced to the classical Lebesgue space
 $L^{p}(\Omega)$ and the Luxemburg norm becomes the standard norm on
$L^{p}(\Omega)$, $$\|\mathbf{u}\|_{L^p(\Omega)}=\left(\int_\Omega |\mathbf{u}(x)|^{p}\,dx\right)^{1/p}.$$

If $p(x)\not\equiv$ constant in $\Omega$, then an important role in manipulating the generalized Lebesgue-Sobolev spaces is played by the modular $\rho_{p(\cdot)}$ of the space $L^{p(\cdot)}(\O)$, which is the mapping $\rho_{p(x)}:L^{p(x)}(\O)\to \R$ defined by
\begin{equation*}
   \rho_{p(x)}(u):=   \int_{\Omega} |u|^{p(x)} \,dx,
\end{equation*}
and the following properties hold:
\begin{gather*}\label{modular1}
|\mathbf{u}|_{p(x)}<1\;\Rightarrow\;|\mathbf{u}|_{p(x)}^{p^+}\leq\rho_{p(x)}(\mathbf{u})\leq|\mathbf{u}|_{p(x)}^{p^-},\\
\label{modular2}
|\mathbf{u}|_{p(x)}>1\;\Rightarrow\;|\mathbf{u}|_{p(x)}^{p^-}\leq\rho_{p(x)}(\mathbf{u})\leq|\mathbf{u}|_{p(x)}^{p^+},\\
\label{modular3}
|\mathbf{u}|_{p(x)}=1\;\Rightarrow\; \rho_{p(x)}(\mathbf{u})=1,\\
\label{modular4}
|\mathbf{u}_n-\mathbf{u}|_{p(x)}\to0\;\Leftrightarrow\;\rho_{p(x)}(\mathbf{u}_n-\mathbf{u})\to 0.
\end{gather*}
For more details about these variable exponent Lebesgue spaces see \cite{bib3,Ham,RD}.

 \begin{rem} \label{rmk} \rm
 Variable exponent Lebesgue spaces resemble the classical Lebesgue spaces in many
 respects, they are separable Banach spaces and the H\"{o}lder inequality holds.
 The inclusions between Lebesgue spaces also naturally generalize, that is,
 if $0<\mbox{meas}(\Omega)<\infty$ and $p(x),q(x)$ are variable exponents such that
 $p(x)<q(x)$ a.e. in $\Omega$, then there exists a continuous embedding
 $L^{q(x)}(\Omega) \hookrightarrow L^{p(x)}(\Omega).$
 \end{rem}

 \begin{definition} \label{def2.3} \rm
 The variable exponent Sobolev space $W^{1,p(x)}$ is defined by
 $$
W^{1,p(x)}(\Omega)=\big\{\mathbf{u}\in L^{p(x)}(\Omega): | \nabla \mathbf{u} | \in
 L^{p(x)}(\Omega) \big\},
$$
 with the norm
\begin{gather*}
\| \mathbf{u}\|_{1,p(x)}=\inf \big\{ \lambda>0: \int_\Omega \Big(\left|
 \frac{\nabla \mathbf{u}}{\lambda}\right|^{p(x)}+\left|\frac{\mathbf{u}}{\lambda}\right|^{p(x)} \Big)\;dx \leq 1 \big\}, \\
\| \mathbf{u} \|_{1,p(x)}=\|\nabla \mathbf{u}\|_{p(x)}+| \mathbf{u} |_{p(x)},
\end{gather*}
 where $$ |\nabla \mathbf{u}|=\sqrt{\sum_{i=1}^N\big(\frac{\partial \mathbf{u}}{\partial x_i}\big)^2}.$$
 \end{definition}

 \begin{definition} \label{def2.8} \rm
 For $p(x)\in C_+(\overline{\Omega})$, let us define the so-called critical Sobolev
 exponent $p^*(x)$ of $p(x)$ by
 $$
 p^*(x)=\begin{cases}
 \frac{3p(x)}{3-p(x)} & \text{if } p(x)<3,\\
 +\infty & \text{if } p(x)\geq 3,
 \end{cases}
 $$
 for every $x\in \overline{\Omega}$.
 \end{definition}
 We define $W^{1,p(x)}_0(\Omega)$ as the closure of $C^{\infty}_0(\Omega)$
with respect to the norm $\|\cdot\|_{p(x)}$,
$$
W^{1,p(x)}_0(\Omega)=\big\{\mathbf{u}: \mathbf{u}|_{\partial\Omega}=0,
 \mathbf{u}\in L^{p(x)}(\Omega), |\nabla \mathbf{u}|\in L^{p(x)}(\Omega)\big\}.
$$
 The dual space of $W_0^{1,p(x)}(\Omega)$ is denoted by $W^{-1,p'(x)}(\Omega)$,
 where
 $$\frac{1}{p(x)}+\frac{1}{p'(x)}=1, \ \
 \mbox{for every} \ \
 x\in \overline{\Omega}.$$

 Next, we recall some embedding results regarding variable exponent
 Lebesgue and Sobolev spaces.

 \begin{theorem}(see \cite[Theorem $1.3$]{FZ1})\label{embedding}
 The following statements hold:
 \begin{itemize}
 \item[(i)] $(W_0^{1,p(x)}(\Omega),\| \cdotp \|)$ is a separable
 and reflexive Banach space.
 \item[(ii)] If $p,q\in C_+(\overline{\Omega})$ and $ q(x)<p^*(x)$ for every $x\in
 \overline{\Omega}$, then there is a compact and continuous embedding
 $$W^{1,p(x)}(\Omega)\hookrightarrow L^{q(x)}(\Omega).$$
 \item[(iii)] There is a constant $C>0$ such that $$|\mathbf{u}|_{p(x)}\leq C|
 |\nabla \mathbf{u}||_{p(x)}, \text{ for all }\mathbf{u}\in W_0^{1,p(x)}(\Omega).$$
 \end{itemize}
 \end{theorem}
Let $$\mathbf{L}^{p(x)}(\O)=L^{p(x)}(\O)\times L^{p(x)}(\O)\times L^{p(x)}(\O)$$ and define
$$\mathbf{W}^{p(x)}(\O)=\{\mathbf{v}\in \mathbf{L}^{p(x)}(\O):\nabla\times \mathbf{v}\in \mathbf{L}^{p(x)}(\O),\nabla\cdot \mathbf{v}=0,\mathbf{v}\cdot \mathbf{n}\mid_{\partial \O}=0\}, $$
where $\mathbf{n}$ denotes the outward unit normal vector to $\partial \O$. Equip $\mathbf{W}^{p(x)}(\O)$ with the norm
$$\|\mathbf{v}\|_{\mathbf{W}^{p(x)}(\O)}=\|\mathbf{v}\|_{L^{p(x)}(\O)}+\|\nabla\times \mathbf{v}\|_{L^{p(x)}(\O)}.$$
If $p^- > 1$, then by Theorem $2.1$ of \cite{Antontsev}, $\mathbf{W}^{p(x)}(\O)$ is a closed subspace of $\mathbf{W}_\mathbf{n}^{1,p(x)}(\O)$, where
$$\mathbf{W}_\mathbf{n}^{1,p(x)}(\O)=\{\mathbf{v}\in W^{1,p(x)}(\O):\mathbf{v}\cdot \mathbf{n}\mid_{\partial \O}=0\}$$
and $$\mathbf{W}^{1,p(x)}(\O)=W^{1,p(x)}(\O)\times W^{1,p(x)}(\O)\times W^{1,p(x)}(\O).$$
Thus, we have the following theorem.
\begin{theorem}\label{ebb1}
(see  \cite[Theorem 2.1]{XWZ}) Assume that $1 < p^- \leq p^+ < \infty$ and $p$ satisfies condition \eqref{cond12}. Then
$\mathbf{W}^{p(x)}(\O)$ is a closed subspace of\; $\mathbf{W}_n^{1,p(x)}(\O)$. Moreover, if $p^- > \frac{6}{5}$, then $\|\nabla\times\cdot\|_{p(x)(\O)}$ is a norm on
$\mathbf{W}^{p(x)}(\O)$
 and there exists $C = C(N, p^-, p^+) > 0$ such that
$$\|\mathbf{v}\|_{W^{1,p(x)}(\O)}\leq C \|\nabla\times \mathbf{v}\|_{L^{p(x)}(\O)}.$$
\end{theorem}
\begin{rem}\label{remm}
By Theorems \ref{embedding} and \ref{ebb1}, the embedding $\mathbf{W}^{p(x)}(\O)\hookrightarrow \mathbf{L}^{q(x)}(\O)$ { is compact}, with $1 < p^- \leq p^+ < 3$, $q\in C(\overline{\O})$, and $1 \leq q(x) < \frac{3p(x)}{3-p(x)}$ in $\overline{\O}$. Moreover,
$(\mathbf{W}^{p(x)}(\O), \|\cdot\|)$ is a uniformly convex, reflexive and separable
Banach space.
\end{rem}
Let
$$\|\mathbf{u}\|_a=\inf\left\{\eta>0:\int_\O\left( \Big|\frac{\nabla\times \mathbf{u}(x)}{\eta}\Big|^{p(x)}+a(x)\Big|\frac{\mathbf{u}(x)}{\eta}\Big|^{p(x)}\right)dx\leq 1\right\}$$
for all $\mathbf{u}\in \mathbf{W}^{p(x)}$. In view of $a^-> 0$ (see condition \ref{hypot:a}), it is easy to see that $\|\cdot\|_a$ is equivalent to the norms
$\|\cdot\|_{\mathbf{W}^{p(x)}(\O)}$ and $\|\cdot\|_{L^{p(x)}(\O)}$ on $\mathbf{W}^{p(x)}(\O)$. In this paper, we shall use for convenience the norm $\|\cdot\|_a$
on the space $\mathbf{W}^{p(x)}(\O)$.
\begin{prop}\label{pip}(see \cite{XWZ})
 Let $$\Lambda_{p(x),a}(\mathbf{u})=\int_\O\left( |\nabla\times \mathbf{u}(x)|^{p(x)}+a(x)|\mathbf{u}(x)|^{p(x)}\right)dx \ \
 \mbox{for all} \ \ \mathbf{u}\in \mathbf{W}^{p(x)}(\O).$$ Then
\begin{enumerate}
  \item $|\mathbf{u}|_{a}<1\;\Rightarrow\;|\mathbf{u}|_{a}^{p^+}\leq\Lambda_{p(x),a}(\mathbf{u})\leq|\mathbf{u}|_{a}^{p^-};$
\item $|\mathbf{u}|_{a}>1\;\Rightarrow\;|\mathbf{u}|_{a}^{p^-}\leq\Lambda_{p(x),a}(\mathbf{u})\leq|\mathbf{u}|_{a}^{p^+}.$
\end{enumerate}
\end{prop}

\begin{prop}\label{pro} The following functional $$\Phi(\mathbf{u})=\int_\O\frac{1}{p(x)}\left( |\nabla\times \mathbf{u}(x)|^{p(x)}+a(x)|\mathbf{u}(x)|^{p(x)}\right)dx,$$ is well defined, even, convex, and sequentially weakly lower semi-continuous. Also, the functional $\Phi$ is of class $C^1$
and $$(\Phi'(\mathbf{u}),\mathbf{v})=\int_\O\left( |\nabla\times \mathbf{u}|^{p(x)-2}\nabla\times \mathbf{u}\cdot\nabla\times \mathbf{v}+a(x)|\mathbf{u}|^{p(x)-2}\mathbf{u}\cdot \mathbf{v}\right)dx,\text{ for all } \mathbf{u}, \mathbf{v}\in\mathbf{W}^{p(x)}(\O),$$
where
 $\langle\cdot,\cdot\rangle$ is the dual pairing between $\mathbf{W}^{p(x)}(\O)$ and its dual $(\mathbf{W}^{p(x)}(\O))^*$.
Similar to \cite{XWZ}, we can deduce that
\begin{enumerate}
  \item [$(i)$] $\Phi': \mathbf{W}^{p(x)}(\O)\to (\mathbf{W}^{p(x)}(\O))^*$ is a continuous, bounded and strictly
monotone operator;
  \item [$(ii)$] $\Phi'$ is a mapping of type $(S_+)$, namely: $\mathbf{u}_n\rightharpoonup \mathbf{u}$ and $\limsup_{n\to\infty} \langle A'(\mathbf{u}_n),\mathbf{u}_n-\mathbf{u}\rangle\leq 0,$
hence $\mathbf{u}_n\to \mathbf{u}$ in $\mathbf{W}^{p(x)}(\O)$.
  \item [$(iii)$] $\Phi': \mathbf{W}^{p(x)}(\O)\to (\mathbf{W}^{p(x)}(\O))^*$ is a homeomorphism.
\end{enumerate}
\end{prop}
\begin{rem} (see \cite[Remark 2.1]{FanD})\label{deriveS+}
We note that the sum of a mapping of type $(S_+)$ and a
weakly-strongly continuous mapping is still a mapping of type $(S_+)$.
Therefore $I'_{\lambda,\mu}=\Phi'-\lambda J'-\mu\Psi'$ is a mapping of type $(S_+)$. Hence any bounded $(P.S)$ sequence of $I_{\lambda,\mu}$ has a convergent
subsequence.
\end{rem}
\subsection{\bf Preliminary lemmas}
From the statement above we know that $\mathbf{W}^{p(x)}(\O)$ is a reflexive and separable
Banach space (see \cite{Antontsev}). Therefore there exist $\{e_j\} \subset \mathbf{W}^{p(x)}(\O)$ and $\{e^*_j\}\subset (\mathbf{W}^{p(x)}(\O))^*$ such that
$$ \mathbf{W}^{p(x)}(\O)=\overline{\mbox{span}}\{e_j:j=1,2,...\},~~~~~~(\mathbf{W}^{p(x)}(\O))^*=\overline{\mbox{span}}\{e^*_j:j=1,2,...\},$$
with \begin{equation*}
      \langle e_j,e_j^*\rangle = \begin{cases}
                                    1, & \mbox{if } i=j \\
                                     0, & \mbox{if } i\neq j.
                                   \end{cases}
     \end{equation*}
     Define
\begin{eqnarray}\label{SOUS}
X_j=\mbox{span}\{e_j\},~~~~~~Y_k=\bigoplus_{j=1}^kX_j,~~~~~~Z_k=\overline{\bigoplus_{j=k}^\infty } X_j.
\end{eqnarray}
We need the following lemmas which will be used in the proof of our main results.
\begin{lemma}\label{BETA}
If $q(x),~~\gamma(x)\in \mathcal{C}_+(\overline{\Omega}),~~q(x),\gamma(x)<p_2^*(x)$ \mbox{ for } $x\in\overline{\Omega},$ let
$$\beta_k=\mbox{sup}\{|\mathbf{u}|_{q(x)}:~\|\mathbf{u}\|_a=1,~\mathbf{u}\in Z_k\},$$
$$\theta_k=\mbox{sup}\{|\mathbf{u}|_{\gamma(x)}:~\|\mathbf{u}\|_a=1,~u\in Z_k\}.$$
Then $lim_{k\rightarrow\infty}\beta_k=0$ and $~~lim_{k\rightarrow\infty}\theta_k=0.$
\end{lemma}
{\bf Proof.} Obviously, $0 < \beta_{k+1} \leq \beta_k$, so $\beta_k \to \beta\geq 0$. Let $\mathbf{u}_k\in Z_k$ satisfy
$$\|\mathbf{u}_k\|_a=1,\quad 0\leq \beta_k-|\mathbf{u}_k|_{q(x)}<\frac1k.$$
Then there exists a subsequence of $\{\mathbf{u}_k\}$ (which we still denote by $\mathbf{u}_k$) such that $\mathbf{u}_k \rightharpoonup \mathbf{u}$, and
$$\langle e^*_j,u\rangle=\lim_{k\to\infty}\langle e^*_j,\mathbf{u}_k\rangle=0,~~\mbox{for all} \ \  e^*_j,$$
which implies that $\mathbf{u}=0$, and so $\mathbf{u}_k \rightharpoonup 0$. Since the embedding from $\mathbf{W}^{p(x)} (\O)$ to $\mathbf{L}^{q(x)} (\O)$ is compact, it follows that $\mathbf{u}_k\to0$ in $\mathbf{L}^{q(x)} (\O)$. Hence, we get $\beta_k\to0$ as $k\to\infty$. The proof for
$\theta_k$ can be obtained by the same procedure.\qed

\begin{lemma}\label{weakly lower semi continuous}
$I_{\lambda,\mu}$  is weakly lower semi-continuous on $\mathbf{W}^{p(x)} (\O)$.
\end{lemma}
\textbf{Proof} By Proposition \ref{pro}, we know that $\Phi$ is weakly lower semi-continuous. Assuming $\mathbf{u}_n\rightharpoonup \mathbf{u}$ in $\mathbf{W}^{p(x)} (\O)$, the compact embedding by Remark \ref{remm} gives us
\begin{eqnarray}\label{h4.2}
\mathbf{u}_n\rightharpoonup~\mathbf{u}~~ \mbox{  in }~~~\mathbf{L}^{p(x)}(\Omega) \mbox{ and } \mathbf{u}_n\rightharpoonup \mathbf{u}~~~~\mbox{  in }~~~~\mathbf{L}^1(\Omega).
\end{eqnarray}
By the mean value theorem, there exists $\mathbf{z}$ which takes on values strictly between $\mathbf{u}$ and $\mathbf{u}_n$ such that
$$\int_{\Omega}|F(x,\mathbf{u}_n)-F(x,\mathbf{u})|dx\leq \int_{\Omega}|\mathbf{u}_n-\mathbf{u}|\sup_{x\in\Omega}|f(x,\mathbf{z})|dx,$$
hence by assumptions \eqref{hypot:F2}, \eqref{hypot:G2}  and \eqref{h4.2}, the functional $J(\mathbf{u}) =\int_{\Omega}F(x,\mathbf{u})dx$ is weakly continuous, and so is $\Psi(\mathbf{u})=\int_{\Omega}G(x,\mathbf{u})dx$. Consequently, the functional $I_{\lambda,\mu}$ is weakly lower semi-continuous.\qed
\section{Proofs of the main results}\label{sect3}
Recall the definitions of $(PS)_c$ and $(PS)^*_c$ conditions.\footnote{We refer the readers to \cite{H1,H2,HH,HHS} for further information on the Palais-Smale condition with assumption weaker than \eqref{hypot:F3}.}
\subsection{\bf The Palais-Smale Compactness Condition}
\begin{definition}
 The $C^1$-functional $I_{\lambda,\mu}$ satisfies the Palais-Smale condition at
 the level $c$ (in short $(PS)_c$ condition) for $c\in \mathbb{R}$ if any sequence
 $(\mathbf{u}_n)_{n\in \mathbb{N}}\subseteq \mathbf{W}^{p(x)}(\O)$ for which $I_{\lambda,\mu}(\mathbf{u}_n)\to c$ and $I'_{\lambda,\mu}(\mathbf{u}_n)\to 0$
 as $n\to \infty$, has a convergent subsequence.
\end{definition}
\begin{definition}
The $C^1$-functional $I_{\lambda,\mu}$ satisfies the $(PS)_c^*$ condition for
 $c\in \mathbb{R}$ if any sequence $(\mathbf{u}_n)_{n\in \mathbb{N}}\subseteq \mathbf{W}^{p(x)}(\O)$ such that $n_j \rightarrow +\infty$, $\mathbf{u}_{n_j} \in Y_{n_j}$, $I_{\lambda,\mu} (\mathbf{u}_{n_j}) \rightarrow c$, and $(I_{\lambda,\mu} |Y_{n_j})'(\mathbf{u}_{n_j} ) \rightarrow 0$,
contains a subsequence converging to a critical point of $I_{\lambda,\mu}$.
\end{definition}

\begin{theorem}\label{PS}
Under the hypotheses of Theorem \ref{th1}, the functional $I_{\lambda,\mu}$ satisfies $(PS)_c$ condition.
\end{theorem}
{\bf Proof. } By Lemma \ref{weakly lower semi continuous} and Remark \ref{deriveS+}, it suffices to verify the boundedness of $(PS)_c$ sequences. Suppose that $(\mathbf{u}_n)_n\subset \mathbf{W}^{p(x)}(\O)$ is a $(PS)$ sequence  at the level $c\in \R$, i.e., $I_{\lambda,\mu}(\mathbf{u}_n)\leq c$
and
$ I'_{\lambda,\mu}(\mathbf{u}_n)\rightarrow 0$ as $n\rightarrow \infty.$ Arguing by contradiction, we assume that $\|\mathbf{u}_n\|_a\to +\infty$. For $n$ large enough, by the conditions \eqref{hypot:F1}, \eqref{hypot:F3}, \eqref{hypot:G1}, \eqref{hypot:G2} and Proposition \ref{pip} we have
\begin{align*}
c+\|\mathbf{u}_n\|_a&\geq  I_{\lambda,\mu}(\mathbf{u}_n)-\frac{1}{\mu}\langle I'_{\lambda,\mu}(\mathbf{u}_n),\mathbf{u}_n\rangle\\&\geq \Bigg(\frac{1}{p^+}-\frac{1}{\mu}\Bigg)\Lambda_{p(x),a}(\mathbf{u}) -\mu\int_{\Omega}\Bigg(\frac{1}{\mu}G(x,\mathbf{u}_n)-g(x,\mathbf{u}_n)\cdot\mathbf{u}_n\Bigg)dx\\&\quad+\lambda\int_{\Omega}\Bigg(\frac{1}{\mu}f(x,\mathbf{u}_n)\cdot\mathbf{u}_n-
F(x,\mathbf{u}_n)\Bigg)dx\\&\geq \Bigg(\frac{1}{p^+}-\frac{1}{\mu}\Bigg)\|\mathbf{u}_n\|^{p^-}_a-\mu\int_{\Omega} \Bigg(\frac{1}{\mu}G(x,\mathbf{u}_n) -g(x,\mathbf{u}_n)\cdot\mathbf{u}_n\Bigg)dx\\&\quad\quad+\lambda\int_{\Omega \cap \{|\mathbf{u}_n|>l\}}\Bigg(\frac{1}{\mu}f(x,\mathbf{u}_n)\cdot\mathbf{u}_n-F(x,\mathbf{u}_n)\Bigg)dx-C|\O|\\
&\geq \Bigg(\frac{1}{p^+}-\frac{1}{\mu}\Bigg)\|\mathbf{u}_n\|^{p^-}_a-C\mu\|\mathbf{u}_n\|^{\gamma^+}_a-C|\O|.
\end{align*}
Dividing the above inequality by $\|\mathbf{u}_n\|^{p^-}_a,$ taking into account that ${p^{-} > \gamma^+}$ and
passing to the limit as $n\to \infty$, we obtain a contradiction. It follows that $(\mathbf{u}_n)_n$ is bounded in $\mathbf{W}^{p(x)}(\O)$. \qed
\begin{theorem}\label{PS*}
Under the hypotheses of Theorem \ref{th2}, functional $I_{\lambda,\mu}$ satisfies $(PS)^*_c$ condition.
\end{theorem}
{\bf Proof.} Suppose that $(\mathbf{u}_{n_j})_j\subset \mathbf{W}^{p(x)}(\O)$ is such that
$$ \mathbf{u}_{n_j}\in Y_{n_j},~~~~I_{\lambda,\mu}(\mathbf{u}_{n_j})\rightarrow c,~~~~~~~(I_{\lambda,\mu}|Y_{n_j})'(\mathbf{u}_{n_j})\rightarrow 0~~~~\mbox{ as }~~n_j\rightarrow +\infty.$$
In a similar way as in the proof of Theorem \ref{PS}, we obtain the boundedness
 of the sequence $(\mathbf{u}_{n_j})_{j\in \mathbb{N^*}}\subseteq \mathbf{W}^{p(x)}(\O)$. Hence, there exists $\mathbf{u}\in \mathbf{W}^{p(x)}(\O)$ such that $\mathbf{u}_{n_j}\rightharpoonup \mathbf{u}$ weakly in $\mathbf{W}^{p(x)}(\O)=\overline{\cup_{n_j}Y_{n_j}}$. Then we can obtain $\mathbf{v}_{n_j}\in Y_{n_j}$ such that $\mathbf{v}_{n_j}\rightharpoonup \mathbf{u}$. We have
$$\langle I'_{\lambda,\mu}(\mathbf{u}_{n_j}), \mathbf{u}_{n_j}-\mathbf{u}\rangle=\langle I'_{\lambda,\mu}(\mathbf{u}_{n_j}), \mathbf{u}_{n_j}-\mathbf{v}_{n_j}\rangle+\langle I'_{\lambda,\mu}(\mathbf{u}_{n_j}), \mathbf{v}_{n_j}-\mathbf{u}\rangle.$$
Since $\mathbf{u}_{n_j}-\mathbf{v}_{n_j}\in Y_{n_j}$, then
$$\langle I'_{\lambda,\mu}(\mathbf{u}_{n_j}), \mathbf{u}_{n_j}-\mathbf{u}\rangle=\langle (I_{\lambda,\mu}|Y_{n_j})'(\mathbf{u}_{n_j}), \mathbf{u}_{n_j}-\mathbf{v}_{n_j}\rangle+\langle I'_{\lambda,\mu}(\mathbf{u}_{n_j}), \mathbf{v}_{n_j}-\mathbf{u}\rangle \rightarrow 0 \mbox{ as } n\rightarrow\infty.$$
Since $I'_{\lambda,\mu}$ is of $S_+$ type, we can deduce that $\mathbf{u}_{n_j}\rightarrow \mathbf{u}$ in $\mathbf{W}^{p(x)}(\O)$. Furthermore, $I'_{\lambda,\mu}(\mathbf{u}_{n_j})\rightarrow I'_{\lambda,\mu}(\mathbf{u})$.
Now we claim that $\mathbf{u}$ is a critical point of $I_{\lambda,\mu}$. Taking $\mathbf{w}_k\in Y_k$, when $n_j\geq k$, we have
$$\langle I'_{\lambda,\mu}(\mathbf{u}),\mathbf{w}_k\rangle =\langle I'_{\lambda,\mu}(\mathbf{u})-I'_{\lambda,\mu}(\mathbf{u}_{n_j}),\mathbf{w}_k\rangle+\langle I'_{\lambda,\mu}(\mathbf{u}_{n_j}),\mathbf{w}_k\rangle$$
$$= \langle I'_{\lambda,\mu}(\mathbf{u})-I'_{\lambda,\mu}(\mathbf{u}_{n_j}),\mathbf{w}_k\rangle+\langle(I_{\lambda,\mu}|Y_{n_j})'(\mathbf{u}_{n_j}),\mathbf{w}_k\rangle.$$
Taking $n_j\rightarrow\infty$, we obtain $ \langle I'_{\lambda,\mu}(\mathbf{u}),\mathbf{w}_k\rangle=0,~~\mbox{for all }\mathbf{w}_k\in Y_k$. So $I'_{\lambda,\mu}(\mathbf{u})=0$, which verifies that $I_{\lambda,\mu}$ satisfies $(PS)^*_c$ condition.\qed

\subsection{\bf Proof of Theorem \ref{th1}}
The following Fountain theorem will be used to get our first result. For the reader convenience, we state it as follows.
\begin{theorem}
\label{Fountain}(Fountain theorem \cite{MW})\label{c}
Let $X$ be a reflexive and separable Banach space, $I\in C^1(X,\mathbb{R})$
 an even functional and let the subspaces $X_k, Y_k, Z_k$ be as defined in
 \eqref{SOUS}.
Suppose that for each $k\in\mathbb{R}$, there exist $\rho_k > r_k > 0$ such that
\begin{enumerate}
\item [$(A_1)$] $\inf\{I(\mathbf{u}):~~
\mathbf{u}\in Z_k,~~||\mathbf{u}||=r\}\rightarrow +\infty ~~\mbox{ as }~~k\rightarrow +\infty;$
\item [$(A_2)$] $\max\{I(\mathbf{u}):~ \mathbf{u}\in Y_k,~~||\mathbf{u}||=\rho_k \}\leq 0;$
\item [$(A_3)$] $I$ satisfies $(PS)$ condition for every $c>0$,
\end{enumerate}
Then $I$ has an unbounded sequence of critical points.
\end{theorem}
{\bf Proof of Theorem \ref{th1}.}  By \eqref{hypot:G4}, \eqref{hypot:F5} and  Theorem \ref{PS}, $I_{\lambda,\mu}$ is an even functional and satisfies  $(PS)_c$ condition. Therefore, by Theorem \ref{Fountain} it suffices to show that if $k$ is large enough, then there exist $\rho_k > r_k > 0$ such that $(A_1)$ and $(A_2)$ hold.
\\{\bf Verification of $(\mathbf{A_1})$:} Let $\mathbf{u}\in Z_k$ with $\|u\|_a>1$.  Then it follows from \eqref{hypot:F2} and  \eqref{hypot:G2} that
\begin{eqnarray*}\label{n1}
I_{\lambda,\mu}(\mathbf{u})&=& \int_\O\frac{1}{p(x)}\left( |\nabla\times \mathbf{u}(x)|^{p(x)}+a(x)|\mathbf{u}(x)|^{p(x)}\right)dx -\lambda\int_\Omega F(x,\mathbf{u})dx-\mu\int_\Omega G(x,\mathbf{u})dx\nonumber\\
&\geq& \frac{1}{p^+}||\mathbf{u}||^{p^-}_a~-\lambda\int_\Omega F(x,\mathbf{u})dx-\mu\int_\Omega G(x,\mathbf{u})dx\nonumber \\
 &\geq& \frac{1}{p^+}||\mathbf{u}||^{p^-}_a-\lambda C\int_\Omega |\mathbf{u}|^{q(x)}dx-\mu C\int_\Omega |\mathbf{u}|^{\gamma(x)}dx\nonumber\\
 &\geq& \frac{1}{p^+}\|\mathbf{u}\|^{p^-}_a-\lambda C\rho_{q(x)}(\mathbf{u})-\mu C\|\mathbf{u}\|^{\gamma^+}_a\\
&\geq&
\begin{cases}
\frac{1}{p^+}\|\mathbf{u}\|^{p^-}_a-C-\mu C\|\mathbf{u}\|^{\gamma^+}_a, & \mbox{if } |\mathbf{u}|_{\gamma(x)}\leq1 \\
\frac{1}{p^+}\|\mathbf{u}\|^{p^-}_a-\frac{\lambda}{q^-} C\beta_k^{q^+}||\mathbf{u}||_a^{q^+}-\mu C\|\mathbf{u}\|^{\gamma^+}_a , & \mbox{if } |\mathbf{u}|_{\gamma(x)}>1.
\end{cases}\\&\geq&\frac{1}{p^+}\|\mathbf{u}\|^{p^-}_a-\frac{\lambda}{q^-} C\beta_k^{q^+}||\mathbf{u}||_a^{q^+}-\mu C\|\mathbf{u}\|^{\gamma^+}_a-C,
\end{eqnarray*}
where \begin{eqnarray*}\label{ww}
 \beta_k=\mbox{sup}\{|\mathbf{u}|_{q(x)}:~\|\mathbf{u}\|_a=1,~\mathbf{u}\in Z_k\}.\end{eqnarray*}
Choose $\|\mathbf{u}\|_a=r_k=\Big(\frac{\lambda}{q^-} q^+ C\beta_k^{q^+}\Big)^{\frac{1}{p^- -q^+}}$
and notice that $p^-<p^+<q^+$.
By Lemma \ref{BETA} we can deduce that $r_k\rightarrow+\infty$ as $k\rightarrow\infty$, hence
\begin{eqnarray*}
I_{\lambda,\mu}(\mathbf{u})&\geq&\frac{1}{p^+}\Big(\frac{\lambda}{q^-} q^+ C\beta_k^{q^+}\Big)^{\frac{p^-}{p^- -q^+}}-\frac{1}{q^+}\Big(\frac{\lambda}{q^-} q^+C\beta_k^{q^+}\Big)\Big(\frac{\lambda}{q^-} q^+ C\beta_k^{q^+}\Big)^{\frac{q^+}{p^- -q^+}}-\mu C\Big(\frac{\lambda}{q^-} q^+ C\beta_k^{q^+}\Big)^{\frac{\gamma^+}{p^- -q^+}}-C.\\&=& \Big(\frac{1}{p^+}-\frac{1}{q^+}\Big)\Big(\frac{\lambda}{q^-} q^+ C\beta_k^{q^+}\Big)^{\frac{p^-}{p^- -q^+}}-\mu C\Big(\frac{\lambda}{q^-} q^+ C\beta_k^{q^+}\Big)^{\frac{\gamma^+}{p^- -q^+}}-C \to \infty \mbox{ as } k\to \infty.
\end{eqnarray*}
\\{\bf Verification of $(\mathbf{A_2})$:}
Clearly, condition \eqref{hypot:F3} implies the existence of two positive constants
$c_1$ and $c_2$ such that
\begin{eqnarray}\label{jh}
F(x,u)\geq c_1|\mathbf{u}|^\theta - c_2,\quad \mbox{for all} \ \ (x,\mathbf{u})\in \O\times \R^3.
\end{eqnarray}
Assume now that \eqref{jh} and \eqref{hypot:G2} hold. Let $\mathbf{u}\in Y_k$ be such that $||\mathbf{u}||_a=\rho_k>r_k>1$. Then
\begin{eqnarray*}
 I_{\lambda,\mu}(\mathbf{u})&=& \int_\O\frac{1}{p(x)}\left( |\nabla\times \mathbf{u}(x)|^{p(x)}+a(x)|\mathbf{u}(x)|^{p(x)}\right)dx -\lambda\int_\Omega F(x,\mathbf{u})dx-\mu\int_\Omega G(x,\mathbf{u})dx\nonumber\\ &\leq& \frac{1}{p^-}\|\mathbf{u}\|^{p^+}_a-\lambda\int_\Omega F(x,\mathbf{u})dx-\mu\int_\Omega G(x,\mathbf{u})dx\nonumber\\
&\leq& \frac{1}{p^-}\|\mathbf{u}\|^{p^+}_a-\lambda c_1\int_{\Omega}|u|^\theta dx+\mu C\int_\Omega|\mathbf{u}|^{\gamma(x)}dx + c_2|\O|.
\end{eqnarray*}
Since $\mbox{dim} Y_k<\infty$, all norms are equivalent in $Y_k$, there are $C^1_W,C^2_W > 0$ such that
\[
\int_{\O}| u| ^{\theta}dx\geq C^1_W\| u\|_a ^{\theta}~\mbox{ and }~\int_{\O}| u| ^{\gamma(x)}dx\leq C^2_W\| u\|_a ^{\gamma^+}.
\]
Hence, we get
$$I_{\lambda,\mu}(\mathbf{u})\leq \frac{1}{p^-}\|\mathbf{u}\|^{p^+}_a -\lambda c_1C^1_W\lambda\|\mathbf{u}\|^\theta_a+\mu CC^2_W\|\mathbf{u}\|^{\gamma^+}_a+c_2|\O|,$$
so we see that $I_{\lambda,\mu}(\mathbf{u})\rightarrow-\infty$ as $\|\mathbf{u}\|_a\rightarrow+\infty$ because $\gamma^+<p^+<\theta$. Conclusion of Theorem \ref{th1} is now reached by invoking Theorem~\ref{c}.\qed
\subsection{\bf Proof of Theorem \ref{th2}}
We shall apply the following Dual Fountain theorem to prove our second main result.
\begin{theorem}
\label{DualFountain} (Dual Fountain theorem  \cite{MW})\label{d}
 Let $X$ be a reflexive and separable Banach space, $I\in C^1(X,\mathbb{R})$
  an even functional, and  $X_k, Y_k, Z_k$ the subspaces defined in \eqref{SOUS}. Assume that there is $k_0 > 0$ such that for each $k > k_0$ , there exist $\rho_k > r_k > 0$ such that
\begin{enumerate}
  \item [$(B_1)$] $\inf\{I(\|\mathbf{u}\|):~\mathbf{u}\in Z_k,~ ||\mathbf{u}||=r_k\}<0;$
  \item [$(B_2)$] $\max\{I(\mathbf{u}):~\mathbf{u}\in Y_k,~ ||\mathbf{u}||=\rho_k\}\geq0;$
  \item [$(B_3)$] $\inf\{I(\mathbf{u}):~\mathbf{u}\in Z_k,~||\mathbf{u}||=\rho_k\}\rightarrow0~~as~~k\rightarrow+\infty;$
  \item [$(B_4)$]  $I$ satisfies $(PS)^*_c$ condition for every $c \in [d_{k_0},0)$.
\end{enumerate}
Then $I$ has a sequence of negative critical values converging to $0$.
\end{theorem}
{\bf Proof of Theorem \ref{th2}.}  According to \eqref{hypot:G4}, \eqref{hypot:F5} and  Theorem \ref{PS*}, $I_{\lambda,\mu}$ is an even functional and satisfies $(PS)^*_c$ condition. Thus it suffices to verify $\mathbf{(B_1)}$, $\mathbf{(B_2)}$ and $\mathbf{(B_3)}$  of Theorem \ref{DualFountain}.
\\\\{\bf Verification of $\mathbf{(B_1)}$:} Assume that \eqref{hypot:F2} and \eqref{hypot:G2} hold. For any $\mathbf{u}\in Z_k$, we have
\begin{eqnarray}
I_{\lambda,\mu}(\mathbf{u}) &=& \int_\O\frac{1}{p(x)}\left( |\nabla\times \mathbf{u}(x)|^{p(x)}+a(x)|\mathbf{u}(x)|^{p(x)}\right)dx -\lambda\int_\Omega F(x,\mathbf{u})dx-\mu\int_\Omega G(x,\mathbf{u})dx\nonumber\\&\geq& \frac{1}{p^+}||\mathbf{u}||^{p^+}_a~-\lambda\int_\Omega F(x,\mathbf{u})dx-\mu\int_\Omega G(x,\mathbf{u})dx\nonumber\\
&\geq& \frac{1}{p^+}||\mathbf{u}||^{p^+}_a-\lambda C||\mathbf{u}||^{q^-}_a-C\mu\int_\Omega|\mathbf{u}|^{\gamma(x)}dx.
\end{eqnarray}
Notice that $q^->p^+$, so there exists small enough $\rho_0>0$ such that $\lambda C||\mathbf{u}||^{q^-}_a\leq \frac{1}{2p^+}||\mathbf{u}||^{p^+}_a$ as $0<\rho=\|\mathbf{u}\|_a\leq\rho_0$.
\\ Then by the proof above, we have
\begin{eqnarray}\label{gh47}
I_{\lambda,\mu}(\mathbf{u})\geq
\begin{cases}
\frac{1}{2p^+}||\mathbf{u}||^{p^+}_a-\mu C\theta_k^{\gamma^-}\|\mathbf{u}\|_a^{\gamma^-}, & \mbox{if } |\mathbf{u}|_{\gamma(x)}\leq1 \\
\frac{1}{2p^+}||\mathbf{u}||^{p^+}_a-\mu C\theta_k^{\gamma^-}\|\mathbf{u}\|_a^{\gamma^+}, & \mbox{if } |\mathbf{u}|_{\gamma(x)}>1.
\end{cases}
\end{eqnarray}
Choose
$$ \rho_k=max\{(2p^+C\mu\theta_k^{\gamma^-})^\frac{1}{p^+-\gamma^-}, (2p^+C\mu\theta_k^{\gamma^+})^\frac{1}{p^+-\gamma^+}\},$$
and notice that $p^+>\gamma^+$, so by Lemma \ref{BETA} we can deduce that $\rho_k\rightarrow0$ as $k\rightarrow\infty$. Hence $I_{\lambda,\mu}(\mathbf{u})\geq0$,\\ i.e., $\mathbf{(B_1)}$ is satisfied.
\\\\{\bf Verification of $\mathbf{(B_2)}$:} Assume that $\mathbf{u}\in Y_k$ with $||\mathbf{u}||_a\leq 1$. Assumption \eqref{hypot:G3} is equivalent to the following
\begin{eqnarray}\label{1abc}
\mbox{there exists } \delta>0,~~G(x,\mathbf{t})\geq C|\mathbf{t}|^\alpha,~~\alpha<p^-, \mbox{ for all }|\mathbf{t}|\in(0,\delta).
\end{eqnarray}
Then by \eqref{1abc} and \eqref{hypot:F2} we have
\begin{eqnarray*}
I_{\lambda,\mu}(\mathbf{u})&=& \int_\O\frac{1}{p(x)}\left( |\nabla\times \mathbf{u}(x)|^{p(x)}+a(x)|\mathbf{u}(x)|^{p(x)}\right)dx -\lambda\int_\Omega F(x,\mathbf{u})dx-\mu\int_\Omega G(x,\mathbf{u})dx\nonumber\\ &\leq & \frac{1}{p^-}||\mathbf{u}||_a^{p^-}+\lambda C\int_\O |\mathbf{u}|^{q(x)}dx-\mu C\int_\O|\mathbf{u}|^\alpha dx\\
&\leq&\frac{1}{p^-}||\mathbf{u}||_a^{p^-}+C||\mathbf{u}||_a^{q^-}-\mu C||\mathbf{u}||_a^\alpha.
\end{eqnarray*}
Since $\alpha<p^-<q^-$, there exists $r_k\in(0,\rho_k)$ such that $I_{\lambda,\mu}(u)<0$ when $\|\mathbf{u}\|_a=r_k$.
\\{\bf Verification of $\mathbf{(B_3)}$:} Notice that $Y_k\cap Z_k\neq\emptyset$ and $r_k<\rho_k$,
so we have\\
$$d_k=\mbox{inf}_{\mathbf{u}\in Z_k,||\mathbf{u}||_a\leq\rho_k}I_{\lambda,\mu}(\mathbf{u})\leq b_k=\mbox{max}_{\mathbf{u}\in Y_k,||\mathbf{u}||_a=r_k} I_{\lambda,\mu}(\mathbf{u})<0.$$
For $\mathbf{u}\in Z_k,~~\|\mathbf{u}\|_a\leq\rho_k$ is small enough. From \eqref{gh47}, we can now obtain
$$I_{\lambda,\mu}(\mathbf{u})\geq\frac{1}{2p^+}||\mathbf{u}||_a^{p^+}-\mu C\theta_{k}^{\gamma^+}||\mathbf{u}||^{\gamma^+}_a.$$
Since $\theta_k\rightarrow 0$ and $k\rightarrow \infty$,
it now follows that $\mathbf{(B_3)}$ is also satisfied.  
 Invoking 
Theorem~\ref{d},
we thus complete the proof of 
Theorem~\ref{th2}.
\qed

\subsection{\bf Proof of Theorem \ref{th3}}
In order to prove Theorem \ref{th3}, we shall need the following two lemmas.
\begin{lemma}\label{lem1.31} For any $\lambda<0,\mu>0$, the following holds:
\begin{enumerate}
  \item [$(1)$] $I_{\lambda,\mu}$ is weakly lower semi-continuous on $\mathbf{W}^{p(x)}(\O)$.
  \item [$(2)$] $I_{\lambda,\mu}$ is bounded from below and coercive on $\mathbf{W}^{p(x)}(\O)$.
\end{enumerate}
\end{lemma}
{\bf Proof.}\\
{\it$(1)$} The proof is similar to that of Lemma \ref{weakly lower semi continuous}, so we shall omit it.
\\
{\it$(2)$} From the hypothese \eqref{hypot:F4}, for any small enough $\epsilon>0$, there exists $M>0$ such that
$$|F(x,\mathbf{t})|\leq \epsilon|\mathbf{t}|^{p^-}~~\mbox{ for } |\mathbf{t}|>M.$$
Therefore, when $\lambda<0,\mu>0$, we can deduce that for any $\mathbf{u}\in \mathbf{W}^{p(x)}(\O)$ with $\|\mathbf{u}\|_a > 1$, the
following holds
\begin{eqnarray}
 I_{\lambda,\mu}(\mathbf{u})&\geq& \frac{1}{p^+}||\mathbf{u}||^{p^-}_a -\int_\Omega F(x,\mathbf{u})dx-\mu \int _\Omega G(x,\mathbf{u})dx\nonumber\\ &\geq& \frac{1}{p^+}||\mathbf{u}||^{p^-}_a+\lambda\epsilon||\mathbf{u}||^{p^-}_a-\mu C||\mathbf{u}||^{\gamma^+}_a.
\end{eqnarray}
Since $\gamma^+<p^-$, $I_{\lambda,\mu}$ is bounded from below and coercive, so $(2)$ is also proved.\qed
\begin{lemma}\label{lem1.32}
Assume that \eqref{hypot:F2} and \eqref{hypot:G3} hold. Then for any
$\lambda<0,\mu>0$ we have
$$\inf_{\mathbf{u}\in \mathbf{W}^{p(x)}(\O)}I_{\lambda,\mu}(\mathbf{u})<0.$$
\end{lemma}
$\\${\bf Proof.}
Using again assumption \eqref{hypot:G3}, there exists
$\delta>0$ such that
\begin{eqnarray}\label{abc}
G(x,\mathbf{e})\geq C\mathbf{e}^\alpha,~~\alpha<p^-,~\  \mbox{for all} \ \  |\mathbf{e}|\in(0,\delta).
\end{eqnarray}
Choose $\mathbf{v}_0\in C^\infty_0(\Omega)$ such that $0<\mathbf{v}_0\leq\delta$, and let $\mathbf{u}_0=s\mathbf{v}_0$. Then by \eqref{hypot:F2} and \eqref{abc}, for $\lambda<0,\mu>0$, we have
 \begin{eqnarray*}
   I_ {\lambda,\mu}(s\mathbf{v}_0) &\leq& s^{p^-}||\mathbf{v}_0||^{p^-}_a-\lambda\int_\Omega F(x,s\mathbf{v}_0)dx-\mu\int_\Omega G(x,s\mathbf{v}_0)dx\\&\leq& s^{p^-}\|\mathbf{v}_0\|^{p^-}_a + |\lambda|C\int_\Omega s^{p(x)}|\mathbf{v}_0|^{p(x)}dx-\mu C\int_\Omega s^\alpha|\mathbf{v}_0|^\alpha dx\nonumber\\&\leq& s^{p^-}||\mathbf{v}_0||^{p^-}_a+|\lambda|C s^{p^-}\int_\Omega|\mathbf{v}_0|^{p(x)}dx-\mu Cs^\alpha\int_\Omega |\mathbf{v}_0|^\alpha dx.
 \end{eqnarray*}
Since $\alpha<p^-$ and $s$ is small enough, it follows that $\inf_{\mathbf{u}\in \mathbf{W}^{p(x)}(\O)}I_{\lambda,\mu}(\mathbf{u})<0$, which completes the proof. \qed
\\
{\bf Proof of Theorem \ref{th3}.} By Lemma \ref{lem1.31}, it follows that for any $\lambda<0,\mu>0$, $I_{\lambda,\mu}$ has a global minimizer $\mathbf{u}_0$ to $I_{\lambda,\mu}(\mathbf{u})$ in $\mathbf{W}^{p(x)}(\O)$  such that $I'_{\lambda,\mu}(\mathbf{u}_0)=0$ (see \cite{MW}). Therefore $\mathbf{u}_0$ is a weak solution of system \eqref{main}.
Moreover, since $I_{\lambda,\mu}(0) = 0$ and $I_{\lambda,\mu}(\mathbf{u}_0) < 0$ (see Lemma \ref{lem1.32}), $\mathbf{u}\neq0$, i.e. $\mathbf{u}_0$ is a nontrivial solution. This completes the proof of Theorem \ref{th3}.\qed
\subsection{\bf Proof of Theorem \ref{th4}}
When $\lambda<0,~\mu<0$, we argue by contradiction that $\mathbf{u}\in \mathbf{W}^{p(x)}(\O)\setminus \{0\}$ is a weak solution of \eqref{main}. Multiplying the first equation  of system \eqref{main} by $\mathbf{u}$, we get
$$ \int_\Omega \nabla\times(|\nabla\times \mathbf{u}|^{p(x)-2}\nabla\times \mathbf{u})\cdot \mathbf{u} dx+\int_\Omega a(x)|\mathbf{u}|^{p(x)-2}\mathbf{u}\cdot \mathbf{u} dx=\lambda\int_\Omega f(x,\mathbf{u})\cdot\mathbf{u}dx+\mu\int_\Omega g(x,\mathbf{u})\cdot\mathbf{u}dx.$$
Using the boundary conditions in \eqref{main} and integrating by parts, we get
$$ \int_\Omega |\nabla\times \mathbf{u}|^{p(x)}dx+\int_\Omega a(x)|\mathbf{u}|^{p(x)}dx=\lambda\int_\Omega f(x,\mathbf{u})\cdot\mathbf{u}dx+\mu\int_\Omega g(x,\mathbf{u})\cdot\mathbf{u}dx,$$
which contradicts  \eqref{hypot:F6} and \eqref{hypot:G5}. This completes the proof of Theorem \ref{th4}.\qed
 
\subsection*{\bf Acknowledgments}
The first author would like to express his
deepest gratitude to the Military School of Aeronautical Specialities, Sfax (ESA) for providing an excellent
atmosphere for  work. The second author was supported by the Slovenian Research Agency grants P1-0292, N1-0114, N1-0083, N1-0064, and J1-0831. The authors wish to acknowledge the referees for several useful comments and valuable suggestions which have helped  improve the presentation.


\begin{thebibliography}{999}
\bibitem{ACX}
 S. Antontsev, M. Chipot, Y. Xie,
\emph{Uniquenesss results for Equation of the  $p(x)$-Laplacian type},
Adv. Math. Sc. Appl.,  $\mathbf{17}$ (1) (2007), 287-304.
\bibitem{AH2}
 T. Adamowicz, P. H\"ast\"o,
\emph{Harnack's inequality and the strong $p(\cdot)$-Laplacian},
J. Differential Equations, $\mathbf{250}$ (2011), 1631-1649.
\bibitem{AMC}
G. A. Afrouzi, M. Mirzapour, N. T. Chung, \emph{Existence and multiplicity of solutions for
Kirchhoff type problems involving $p(x)$-biharmonic operators}, Z. Anal. Anwend., $\mathbf{33}$ (2014), 289-303.
\bibitem{AMS0}
S. Antontsev, F. Miranda, L. Santos, \emph{A class of electromagnetic $p$-curl systems: blow-up and finite time extinction}, Nonlinear Analysis: Theory, Methods and Applications, $\mathbf{75}$(9) (2012), 3916-3929.
\bibitem{AMS}
 R. Aboulaich, D. Meskine, A. Souissi, \emph{New diffusion models in image processing},
Comput. Math. Appl., $\mathbf{56}$ (2008), 874-882.
\bibitem{Antontsev}
S. Antontsev, F. Miranda, L. Santos, \emph{Blow-up and finite time extinction for $p(x,t)$-curl systems
arising in electromagnetism}, J. Math. Anal. Appl., $\mathbf{440}$ (2016), 300-322.
\bibitem{AR} S. N. Antontsev, J. F. Rodrigues,
\emph{On stationary thermorheological viscous flows},
Ann. Univ. Ferrara Sez. VII Sci. Mat., $\mathbf{52}$ (1) (2006), 19-36.
\bibitem{AS}
S. N. Antontsev,  S. I. Shmarev, \emph{A model porous medium equation with variable exponent of nonlinearity: Existence,
uniqueness and localization properties of solutions}, Nonlinear Anal., $\mathbf{60}$ (2005), 515-545.
\bibitem{BD}
A. Bahrouni, D. Repov\v{s}, \emph{Existence and nonexistence of solutions for $p(x)$-curl systems arising in electromagnetism}, Complex Variables and Elliptic Equations, $\mathbf{63}$ (2018), 292-301.
\bibitem{BJF}
B. Ge, J. F. Lu, \emph{Existence and Multiplicity of Solutions for $p (x)$-Curl Systems Without the Ambrosetti-Rabinowitz Condition}, Mediterranean Journal of Mathematics, $\mathbf{16}$(2) (2019), p.45
\bibitem{BO}
 S. Byun, J. Ok,
\emph{On $W^{1,q(x)}$-estimates for elliptic equations of $p(x)$-Laplacian type},
J. Math. Pures Appl., $\mathbf{106}$ (3) (2016), 512-545.
\bibitem{Chung1}
N. T. Chung, \emph{ Existence of solutions for perturbed fourth order elliptic equations with variable exponents},
Electron. J. Qual. Theory Differ. Equ., $\mathbf{2018}$ (96) (2018), 1-19.
\bibitem{Chung2}
N. T. Chung, \emph{ Some remarks on a class of $p(x)$-Laplacian Robin eigenvalue problems}, Mediterr. J. Math.,
$\mathbf{15}$(4) (2018), p.147
\bibitem{NTC1}
N. T. Chung, \emph{Multiple solutions for a $p(x)$-Kirchhoff-type equation with sign-changing
nonlinearities}, Complex Var. Elliptic Equ., $\mathbf{58}$(12) (2013), 1637-1646.
\bibitem{JF}
 J. Chabrowski, Y. Fu, \emph{Existence of solutions for $p(x)$-Laplacian problems on a bounded domain},
J. Math. Anal. Appl., $\mathbf{306}$ (2005), 604-618.
\bibitem{CLR}  Y. Chen, S. Levine, M. Rao,
\emph{Variable exponent, linear growth functionals in image restoration},
SIAM J. Appl. Math., $\mathbf{66}$(4) (2006), 1383-1406.
\bibitem{bib3}
D. E. Edmunds, J. Lang, A. Nekvinda, \emph{On $L^{p(x)}$ norms}, Proc. Roy. Soc. London Ser. A, \textbf{455} (1999),
219-225.
\bibitem{FZ1}
 X. L. Fan, D. Zhao, \emph{On the spaces $L^{p(x)}$ and $W^{m,p(x)}$}, J. Math. Anal. Appl., $\mathbf{263}$ (2001), 424-446.
\bibitem{FZ}
 Y. Q. Fu, X. Zhang,
\emph{Multiple solutions for a class of $p(x)$-Laplacian equations in
$\mathbb{R}^N$ involving the critical exponent},
 Proc. R. Soc.  A,  $\mathbf{466}$ (2010), 1667-1686.
\bibitem{FanD}
X. L. Fan, S. G. Deng, \emph{Remarks on Ricceri's variational principle and applications to the $p(x)-$Laplacian equations},
Nonlinear Anal., $\mathbf{67}$ (2007), 3064-3075.
\bibitem{FanZ} X. L. Fan, D. Zhao,
\emph{Nodal solutions of $p(x)$-Laplacian equations},
Nonlinear Anal., {\bf 67} (2007), 2859--2868.
\bibitem{GL}
 B.  Ge, L. L. Liu,
\emph{Infinitely many solutions for differential inclusion problems in
$\mathbb R^N$ involving the $p(x)$-Laplacian},
Z. Angew. Math. Phys., $\mathbf{67}$ (1) (2016), p.16
\bibitem{Ham}
M. K. Hamdani, A. Harrabi, F. Mtiri, and D. D. Repov\v{s}, \emph{Existence and multiplicity results for a new $p(x)$-Kirchhoff problem}.  Nonlinear Analysis $\mathbf{190}$ (2020): 111598.
\bibitem{H1}
M.~K. Hamdani, \emph{Multiple solutions for Grushin operator without odd nonlinearity}. Asian-European J. Math., (2019),
doi:10.1142/S1793557120501314
\bibitem{H2}
M. K. Hamdani, \emph{On nonlocal asymmetric Kirchhoff problem,} Asian-European J. Math., (2019), doi: 10.1142/S1793557120300018
\bibitem{HH}
M. K. Hamdani, A. Harrabi, \emph{High-order Kirchhoff problems in bounded and unbounded domains,} (arXiv:1807.11040v3 [math.AP]), 5 Aug 2019.
\bibitem{Harrabi}
M. K. Hamdani, A. Harrabi, \emph{ $L^\infty$-norm Estimates of Weak Solutions via their Morse indices for the $m$-Laplacian Problems}. Results in Mathematics, $\mathbf{74}$(1) (2019), p. 69
\bibitem{HHS}
A. Harrabi, M. K. Hamdani, A. Selmi, \emph{Existence results of the zero mass polyharmonic system,} Complex Variables and Elliptic Equations (2019), doi:10.1080/17476933.2019.1679794
\bibitem{HM}
Y. Hong-Ming, \emph{Regularity of weak solution to a $ p $-curl-system}, Differential and Integral Equations, $\mathbf{19}$(4) (2006), 361-368.
\bibitem{KR} O. Kov\'a\v{c}ik, J. R\'akosn\'ik;
\emph{On spaces $L^{p(x)}$ and $W^{k,\,p(x)}$},
Czechoslovak Math. J., $\mathbf{41}$ (1991), 592-618.
\bibitem{LLP}
F. Li, Z. Li, L. Pi,
\emph{Variable exponent functionals in image restoration},
 Appl. Math. Comput., $\mathbf{216}$ (3) (2010), 870--882.
\bibitem{MR}
M. Mih\v{a}silescu, V. D.  R\v{a}dulescu,
\emph{A multiplicity result for a nonlinear degenerate problem arising in
the theory of electrorheological fluids}, Proc. R. Soc. A $\mathbf{462}$ (2006),
2625-2641.
\bibitem{MR2007}
M. Mih\v{a}ilescu, V. D.  R\v{a}dulescu,
\emph{On a nonhomogeneous quasilinear eigenvalue problem in Sobolev spaces with
variable exponent}, Proc. Amer. Math. Soc., $\mathbf{135}$ (9) (2007), 2929-2937.
\bibitem{RR1996} K. R. Rajagopal, M. R\.{u}\u{z}i\u{c}ka,
\emph{On the modeling of electrorheological materials},
Mech. Research Comm., $\mathbf{23}$ (1996) 401-407.
\bibitem{RR2001}K. R. Rajagopal, M. R\.{u}\u{z}i\u{c}ka,
\emph{Mathematical modeling of electrorheological materials},
Cont. Mech. and Thermodynamics, $\mathbf{13}$ (2001), 59-78.
\bibitem{Dusan} D. D. Repov\v{s},  \emph{Stationary waves of Schr\"{o}dinger-type equations with variable exponent}, Anal. Appl. (Singap.), $\mathbf{13}$ (2015), 645-661.
\bibitem{RD}
V. D. R\u{a}dulescu, D. D. Repov\v{s}, \emph{Partial differential equations with variable exponents: variational
methods and qualitative analysis}, CRC Press, Boca Raton, 2015.
\bibitem{R} M. R\.{u}\u{z}i\u{c}ka, \emph{Electro-Rheological Fluids: Modeling and
Mathematical Theory}, Springer-Verlag Berlin, 2000.
\bibitem{XWZ}
M. Q. Xiang, F. L. Wang,  B. L. Zhang, \emph{Existence and multiplicity of solutions
for $p(x)$-curl systems arising in electromagnetism}, Journal of Mathematical Analysis and Applications, $\mathbf{448}$ (2017), 1600-1617.
\bibitem{MW}
 M. Willem, \emph{Minimax Theorems}, Birkh\"{a}user, Boston, 1996.
\end{thebibliography}
\end{document}